\documentclass{ifacconf}

\usepackage{graphicx}      
\usepackage{amsmath,amssymb}
\usepackage{optidef}
\usepackage{xcolor, amsmath, steinmetz, float,booktabs,subcaption}
\usepackage[width=.9\textwidth]{caption}
\usepackage[utf8x]{inputenc}
\usepackage[british,UKenglish,USenglish,american]{babel}
\usepackage{natbib}     
\usepackage{algorithm}
\usepackage{algpseudocode}
\usepackage{pgfplots}
\usepackage{grffile}
\pgfplotsset{compat=newest}
\usetikzlibrary{plotmarks}
\usetikzlibrary{arrows.meta}
\usepgfplotslibrary{patchplots}
\usepackage{tikz}
\usetikzlibrary{shapes.geometric}
\usetikzlibrary{shapes.arrows}
\usepackage{array}
\usepackage{amssymb}

\newtheorem{assumption}{Assumption}
\DeclareMathOperator*{\argmin}{argmin}
\begin{document}
 \renewcommand{\abstractname}{\textbf{Abstract: }}
\definecolor{mycolor1}{rgb}{0.00000,0.44700,0.74100}%
\definecolor{mycolor2}{rgb}{0.85000,0.32500,0.09800}%
\definecolor{mycolor3}{rgb}{0.92900,0.69400,0.12500}%
\definecolor{mycolor4}{rgb}{0.49400,0.18400,0.55600}%
\definecolor{mycolor5}{rgb}{0.46600,0.67400,0.18800}%
\definecolor{mycolor6}{rgb}{0.30100,0.74500,0.93300}%
\definecolor{mycolor7}{rgb}{0.63500,0.07800,0.18400}%
\begin{frontmatter}

\title{A Method for Distributed Transactive Control in Power Systems based on the Projected Consensus Algorithm} 

\author[First]{Eder Baron-Prada} 
\author[Second]{C\'esar A. Uribe}
\author[First]{Eduardo Mojica-Nava}

\address[First]{Universidad Nacional de Colombia,\\ 
   e-mail: \{edbaronp,eamojican\}@unal.edu.co.}
\address[Second]{ECE Department and Coordinated Science Lab \\University of Illinois at Urbana-Champaign,\\ 
   e-mail: cauribe2@illinois.edu.}

\begin{abstract}                 
The shift of power systems toward a smarter grid has brought devices such as distributed generators and smart loads with an increase of the operational challenges for the system operator. These challenges are related to the real-time implementation as well as control and stability issues. We present a distributed transactive control strategy, based on the projected consensus algorithm, to operate the distributed energy resources and smart loads of a power system toward optimal social welfare. We consider two types of agents: Generators, and smart loads. Each agent iteratively optimizes its local utility function based on local information obtained from its neighbors and global information obtained through the network of agents. We show convergence analysis and numerical results for the proposed method. 
\end{abstract}


\end{frontmatter}
\section{Introduction} 

The increased use of devices for monitoring and controlling power networks poses several coordination challenges to assure a robust operation~\citep{Kok2016,Cherukuri2016}. Thus, control theory and optimization algorithms have become one of the main tools to deal with these challenges \citep{Hale2017,Nedic2015,uribe2017a,uribe2017b}. Centralized approaches often suffer from computation and communication overheads. These communication requirements are evident in large-scale systems. Distributed control approaches eases the communication costs of large-scale systems. Several approaches in distributed control for economic dispatch problems have been used to achieve optimal power flow and voltage control \citep{Cherukuri2017,Mojica-Nava2016}. However, these algorithms are subject to several technical requirements.  For example, the coordination 
of devices/agents (e.g., generators and loads) is required to maintain the stability of the power system, i.e., maintain equilibria between power generated and power demanded.

Recently, some strategies has been proposed for design decentralized feedback controllers that steer the system to the optimal solution without explicitly solving the economic dispatch problem \citep{Li2017}, there are also some approaches that consider a more complete power flow problem in a distributed fashion with communications constraints and losses of control signals \citep{Simonetto2016}. Furthermore, several demand-response strategies have been proposed to solve the economic dispatch problem with changes in the load and the generation \citep{Knudsen2016,Shiltz2016,Bejestani2014}. In this context, transactive control has been shown effective to assure coordination of a vast number of devices, including smart loads \citep{Kok2016}. Transactive control uses a market mechanism that allows agents to interact through an economic signal to properly distribute the available resources \citep{Kok2016}.

In this paper, we propose a distributed transactive control method based on consensus-based constrained optimization~\citep{Ozdaglar2010}. Particularly, we use the distributed Minimum Diameter Spanning Tree (MDST) algorithm \citep{Bui2004} to share some required global network parameters in a distributed way (e.g., equality between the load demanded and power generated). Then, we use the projected distributed consensus method to optimize the social welfare. This approach does not need confidential information to be shared, such as the gradient of local utility functions or incremental costs. This allows the implementation of the proposed algorithm on power systems with various proprietary network infrastructure agents. An example of these systems are distribution systems with microgrids where a distribution system usually has an independent system operator (ISO) as owner and microgrids are users' property. Contrary to recent literature~\citep{Cortes2016,Mojica-Nava2014}, the main contribution of this paper is to remove the assumption that all agents have immediate and global access to the {power demand, the power generation and the number of agents} beforehand. Our proposed method allows all agents to obtain these global parameters and achieve the correct operation of the power system, minimize the local cost of each agent and satisfy the local constraints.

The rest of the paper is organized as follows. Section~\ref{sec:problem} presents the problem statement, where we formulate the distributed transactive optimization problem. Section~\ref{sec:DED} describes the proposed transactive control algorithm. Section~\ref{sec:cases} presents experimental results to test the behavior of the proposed algorithm. Finally, Section~\ref{sec:conclusions} shows conclusions and future work.    

\section{Problem Formulation}\label{sec:problem}

We consider a transactive grid, where there are two classes of agents, generators and consumers, seeking to optimize its utility function. The network of agents is represented as a set of nodes in a graph, where edges indicate the communication links among agents. 

\subsection{Preliminaries: Graph Theory}\label{sec:graph}

Let $\mathcal{G}=(\mathcal{V},\mathcal{E})$ be a connected, undirected and unweighted graph \cite[Section 3.4]{Bullo2018}, where $\mathcal{V}=\{1,2,3,...,M+N\} $, and $M$ and $N$ are the number of generators and loads respectively. In addition, $\mathcal{E}$ is the set of communication links between agents, i.e., $(a,b) \in \mathcal{E}$ if there is a link between $a \in \mathcal{V}$ and $b \in \mathcal{V}$. Besides, the neighbors of an agent $a \in \mathcal{V}$ are denoted by $\mathcal{N}_a=\{ b | (a,b) \in \mathcal{E}\}$. The graph $\mathcal{G}$ is connected, which means that there exists at least one path between any two distinct nodes \cite[Section 2.2]{west2000}. Additionally, connectedness implies the existence of at least one spanning tree \cite[subsection 3.2]{Bullo2018}. 

\subsection{Generators and Consumers Model}

Generators are considered as distributed energy resources. We assume that all generators are dispatchable. A dispatchable generator can change the power that it is generating by taking into account the system requirements (i.e., the generation-demand balance), and respond dynamically to changes in the power demanded by consumers. The set of generators is denoted as $G=\{1,2,3,...,N\}$. Moreover, the cost of generation is assumed quadratic \citep{Shiltz2016}, i.e.,     

\begin{align*}
C(P_{g_i})= \rho_{g_i} P_{g_i}+\frac{\beta_{g_i}}{2} P_{g_i}^2,
\end{align*}

where $\rho_g\in \mathbb{R}^{N}$ and $\beta_g\in \mathbb{R}^{N}$ are cost coefficients, $P_{g}(k)\in \mathbb{R}^{N}$ is a vector that contains the power delivered by each generator at time instant $k$, $P_{g}(k)=[ P_{g_1}(k), P_{g_2}(k), ..., P_{g_N}(k)]^{\top}$. Moreover, each generator has local constraints given by

\begin{subequations}
        \begin{align}
        \overline{P_{g_i}} \geq P_{g_i}\geq\underline{P_{g_i}},
        \label{eqn:gen_constraints}
        \end{align}
        \begin{align}
        \overline{T_{g_{i}}}\geq P_{g_i}(k) - P_{g_i}(k-1)\geq\underline{T_{g_{i}}},
        \label{eqn:Tc_gen_constraints}
        \end{align}
    \end{subequations}
where $\underline{P_{g_i}}$ and $\overline{P_{g_i}}$ are the minimum and maximum power delivered by the $i$-th generator, respectively. Additionally, $\underline{T_{g_{i}}}$ and $\overline{T_{g_{i}}}$ are the rate constraints for the $i$-th generator. Rate constraints impose a limit on the rate of change in power generated due to physical limitations. In order to include the power cost of each generator and maximum power capacity we rewrite the cost function as 

\begin{align}
    C(P_{g_i})=  \frac{1}{c_i}\left(P_{g_i}-\frac{P_{g_i}^2}{\overline{P_{g_i}}}  \right),
    \label{eqn:costgen}
    \end{align}
 where $\beta_{g_i} ={-1}/{\overline{P_{g_i}}c_i}$, $\rho_{g_i}={1}/{c_i}$, and ${c_i}$ is the power cost. The set of consumers is denoted as $ {D}=\{1,2,3,...,M\}$. Consumers are assumed to be agents that can obtain the system power cost through demand response devices. These devices are capable of making decisions about the amount of power demanded by each consumer with the objective of maximizing its utility. Consumers are assumed to have controllable loads. The base load is the amount of power that each user consume and can change in time without any explanation. The power consumption that is adapted according to the system parameters is called the variable load. The consumers utility function is defined as
 
 \begin{align*}
 U(P_{d_j})= \rho_{d_i} P_{d_j}+\frac{\beta_{d_i}}{2} P_{d_j}^2,
 \end{align*}
 
where $\rho_d\in \mathbb{R}^{M}$ and $\beta_d\in \mathbb{R}^{M}$ are utility coefficients, $P_{d}(k)\in \mathbb{R}^{M}$ is the vector that contains the load demanded by each smart consumer at time instant $k$, $P_{d}(k)=[ P_{d_1}(k), P_{d_2}(k), ..., P_{d_M}(k)]^\top$.  The agent $j$ has local constraints given by

\begin{subequations}
		\begin{align}
		 \overline{P_{d_j}} \geq P_{d_j}\geq \underline{P_{d_j}}, \quad
		\label{eqn:con_constraints}
		\end{align}
		\begin{align}
		\overline{T_{d_{j}}}\geq P_{d_j}(k) - P_{d_j}(k-1)\geq \underline{T_{d_{j}}},
		\label{eqn:Tc_con_constraints}
		\end{align}
	\end{subequations}
where $\underline{P_{d_j}}$ and $\overline{P_{d_j}}$ are the minimum and maximum power demanded by the load $j$-th. $\underline{T_{d_{j}}}$ and $\overline{T_{d_{j}}}$ are the rate constraints for the $j$-th load, these constraints limit the rate of load change due to physical restrictions in power consumption devices. To include the power cost in the system and the maximum power load between the utility function, we transform the cost function as

\begin{align}
    U(P_{d_j})=  \frac{V_{d_j}}{c_\Omega}\left(P_{d_j}-\frac{P_{d_j}^2}{\overline{P_{d_j}}}  \right),
    \label{eqn:uticonsumers}
    \end{align}
    where $\beta_{d_j} ={-V_{d_j}}/{\overline{P_{d_j}}}$, $\rho_{d_j}={V_{d_j}}/{c_\Omega}$, and $V_{d_j}$ is the value of the power for the agent $j$. These values can change as a function of the power consumption preferences of each consumer. $c_\Omega$ is a global variable that indicates the power cost in the system and it is defined as $c_\Omega={\sum_{i=0}^{N} P_{g_i}(k)c_i}/{\sum_{j=1}^M P_{d_j}(k)}$, where $c_i$ is the power cost of agent $i$. 

Note that the utility functions in \eqref{eqn:costgen} and \eqref{eqn:uticonsumers} are chosen as the integral of logistic-type function. If the power consumed is not at the set point, the gradient obtained from the utility function will have a larger magnitude.  However, if the power consumed is close to the set point, the gradient will have a lower magnitude \citep{Britton2003}.

\subsection{Social Welfare Optimization Problem} \label{subsec:SWO}

The social welfare problem is related to optimizing the state for each agent in the power system, i.e., to maximize the utility of the consumers and to minimize the cost of the generators, and it is obtained from \eqref{eqn:costgen} and \eqref{eqn:uticonsumers} as

\begin{align}
    S_W&=\sum_{j=1}^{M}U(P_{d_j})-\sum_{i=1}^{N}C(P_{g_i}).
    \label{eqn:SW}
\end{align}
  
The main goal is maximize \eqref{eqn:SW} as it is shown in \eqref{eqn:opt}, taking into account constraints \eqref{eqn:gen_constraints}, \eqref{eqn:con_constraints}. 

\begin{maxi!}[2]
        {P_{d_j}, P_{g_i} {i\in G, j\in D}}{S_W}     
        {\label{eqn:opt}}{}
        \addConstraint{\sum_{j=1}^M P_{d_j}-\sum_{i=0}^{N} P_{g_i}}{=0} \label{eqn:global_cons}
        \addConstraint{\overline{P_{g_i}}\geq P_{g_i} \geq \underline{ P_{g_i}}}{\quad \forall i  \in G} 
        \addConstraint{\overline{P_{d_j}}\geq P_{d_j} \geq \underline{ P_{d_j}}}{\quad \forall j  \in D}. 
        \end{maxi!}

{Constraints \eqref{eqn:Tc_gen_constraints} and \eqref{eqn:Tc_con_constraints} are addressed with the stepsize of the Algorithm \ref{alg:distributed_transactive_control} explained in subsection \ref{subsec:convergence}.} We assume that the local objective function and the local constraint set are known to an agent only. Problem \eqref{eqn:opt} is strongly convex and local constraints are closed convex sets. Finally, it is assumed that the feasible set is non-empty. In the following section, we present the preliminary concepts for solving the optimization problem proposed in this section in a distributed way.

	\section{Distributed Transactive Control}\label{sec:DED}
	We propose a distributed transactive control method based on the projected distributed gradient descent method \citep{Ozdaglar2010}. In order to find a solution the following assumptions are made. These assumptions are related to graph properties, communications and agents initial knowledge.
	
 \begin{assumption}
The set of feasible points for problem the \eqref{eqn:opt} is non-empty.
\label{as:op_set_ne}
\end{assumption}   

\begin{assumption}
	The graph $\mathcal{G}$ is connected, static, undirected, and unweighted. Links are assumed lossless, without delays and synchronous.\label{as:graph}
\end{assumption}

\begin{assumption}
The agents know their set of neighbors and hence their own cardinality, $r_a$, where $r_a= |\mathcal{N}_a|$.
\label{as:neighbors}
\end{assumption}

To solve \eqref{eqn:opt}, it is necessary that all agents have information about the global state of the network. Agents need the following data: total power demanded ($P_D$), total power delivered ($P_G$), power cost in the system ($c_\Omega$) and finally the number of generators and consumers in the network ($N$) and ($M$) respectively. Where ($P_D$) is calculated such as $P_D=\sum_{j=1}^M P_{d_j}(k)$ and ($P_G$) is calculated such as $P_G=\sum_{i=1}^N P_{g_i}(k)$. Previous results in the literature assume that global parameters can be obtained \textit{a priori} by agents. This implies that a central entity obtains the global parameters and sends them to all the agents in the network agents as in \cite{Cortes2016,Mojica-Nava2014}. To avoid the requirement of the central entity and construct a fully distributed algorithm, we use an efficient algorithm to find the global parameters in a distributed way. We use the {Minimum Diameter Spanning Tree} (MDST) algorithm to reach consensus in a finite number of steps, those algorithms are explained in subsection \ref{subsec:DAMDST} and \ref{subsec:FTDA} respectively.

The MDST algorithm is used to share the global parameters with all the agents in the system by using a spanning tree of the graph $\mathcal{G}$. Later, the same spanning tree is used to execute the distributed projection gradient algorithm to compute the iterations of the optimization variables. 

\subsection{Distributed Algorithm for the MDST} \label{subsec:DAMDST}

We use the finite time consensus algorithm to achieve common knowledge of the global system's parameters at every node \citep{Mou2014}. However, this algorithm only works on spanning tree graphs. Finding a spanning tree in a graph is a problem heavily studied in recent years \citep{Elkin2006,Gfeller2011,Bui2004}. We use the approach in \cite{Bui2004} to guarantee the best convergence time for the finite-time distributed algorithm, because this algorithm converges in maximum $d$ steps, being $d$ the graph diameter. 

    \begin{lem}\label{lemma:MDST}
        \cite[Theorem 5]{Bui2004}: Consider a graph $\mathcal{G}$ and let Assumptions \ref{as:graph} and \ref{as:neighbors} hold on $\mathcal{G}$. Then, the distributed algorithm for the MDST proposed in \cite[Theorem 5]{Bui2004} finds a MDST of $\mathcal{G}$, in $O(n)$ iterations.
    \end{lem}

\begin{pf}
The distributed algorithm for the MDST has to calculate the All-Pairs Shortest Path (APSP) of the network. The time of execution of APSP is $O(n)$ \cite[Lemma 4]{Bui2004}. Then, each node knows which node is the shortest path to another node. After that, the absolute center of the graph is calculated, given a node with the lowest {eccentricity} $a_{min}$, the information about which node is the center of the graph is sent to another nodes in at most $O(n)$. Now consider the collection of all paths produced by APSP that begins in any node in the network and end in $a_{min}$, the set of path forms a tree rooted in $a_{min}$, which is the MDST of $\mathcal{G}$. Therefore each node knows a route to $a_{min}$ and the MDST is built through of knowledge of the shortest path to $a_{min}$ for all nodes~\cite[Subsection 2.1.5]{Bui2004}. $\blacksquare$
     \end{pf}

Once the MDST is created the global parameters are calculated through the distributed algorithm for the MDST as described in the following subsection.

\subsection{Finite-time Distributed Averaging} \label{subsec:FTDA}

The algorithm proposed in \cite{Mou2014} is used to calculate the global parameters $P_D$, $P_G$, $N$ and $M$. This algorithm can achieve consensus in finite time for spanning trees graphs. Initially, each agent has initial values for the global parameters as follows:

\begin{subequations} \label{eqn:GV_ini}
		\begin{align}\label{eqn:PG_ini}
		P_{G_a}(0)= 
		\begin{cases}
		P_{g_{a}}(0),& \text{if } a {\in G},\\
		0,              & \text{otherwise},
		\end{cases}
		\end{align}
		\begin{align}\label{eqn:PD_ini}
		P_{D_a}(0)= 
		\begin{cases}
		P_{d_{a}}(0),& \text{if } a {\in D},\\
		0,              & \text{otherwise},
		\end{cases}
		\end{align}
		\begin{align}\label{eqn:N_ini}
		N_a(0)= 
		\begin{cases}
		1,& \text{if } a {\in G},\\
		0,              & \text{otherwise},
		\end{cases}
		\end{align}
        \begin{align}\label{eqn:M_ini}
		M_a(0)= 
		\begin{cases}
		1,& \text{if } a {\in D},\\
		0,              & \text{otherwise},
		\end{cases}
		\end{align}
\end{subequations}

where $N_a$ and $M_a$ are the number of generators and consumers respectively known to the agent $a$ before starting the algorithm. Let $x_a$ be any global parameter previously presented in \eqref{eqn:GV_ini}, for each of previous values each agent, performs in parallel the following update action:

\begin{align}\label{eqn:update_GV}
		\begin{split}
				x_a(q+1)\hspace{6cm}\\= 
			\begin{cases}
				x_a(0)+ \sum_{b\in \mathcal{N}_{a}}x_b(0),			& \text{if } q=0;\\
				\sum_{b\in \mathcal{N}_{a}}x_b(q)+(1-r_a)x_a(q-1),  & \text{if } q \geq 1.
			\end{cases}
		\end{split}
	\end{align}
    
   Variable $q$ is used to represent steps. Algorithm \ref{alg:global-average} explains how \eqref{eqn:update_GV} is used to calculate the global parameters. 

 \begin{algorithm}
 \begin{algorithmic}[1]
 	\State \textbf{Executed by:}  Agents $a \in \mathcal{V}=\{1,...,N+M\}$
 	\Require Spanning Tree Neighbors
    \State \textbf{Initialize:} $P_{D_a}(0)$, $P_{G_a}(0)$, $N_{a}(0)$, $M_{a}(0)$ and Set $q=0$
    \State $P_{D_a}(1)=P_{D_a}(0)+ \sum_{b\in \mathcal{N}_{a}}P_{D_b}(0)$
                \State $P_{G_a}(1)=P_{G_a}(0)+ \sum_{b\in \mathcal{N}_{a}}P_{G_b}(0)$
                \State $N_{a}(1)=N_{a}(0)+ \sum_{b\in \mathcal{N}_{a}}N_{b}(0)$
                \State $M_{a}(1)=M_{a}(0)+ \sum_{b\in \mathcal{M}_{a}}M_{b}(0)$
 		\While{$q \leq d$}
 		\State Send $P_{D_a}(q)$, $P_{G_a}(q)$,$N_{a}(q)$ and $M_{a}(q)$ to Spanning Tree neighbors
             	\State $P_{D_a}(q+1)=\sum_{b\in \mathcal{N}_{a}}P_{D_b}(q)+(1-r_a)P_{D_a}(q-1)$
                \State $P_{G_a}(q+1)=\sum_{b\in \mathcal{N}_{a}}P_{G_b}(q)+(1-r_a)P_{G_a}(q-1)$
                \State $N_{a}(q+1)=\sum_{b\in \mathcal{N}_{a}}N_{b}(q)+(1-r_a)N_{a}(q-1)$
                \State $M_{a}(q+1)=\sum_{b\in \mathcal{N}_{a}}M_{b}(q)+(1-r_a)M_{a}(q-1)$
        \State Set $q=q+1$;
 		\EndWhile
 \State \textbf{return} $P_{D_a}(q)$, $P_{G_a}(q)$, $N_a(q)$  and $M_a(q)$ 
 \State \Comment{ $P_{D_a}(q)$=$P_D$, $P_{G_a}(q)$=$P_G$, $N_{a}(q)$=$N$, $M_{a}(q)$=$M$}
 \end{algorithmic}
\caption{Finite-time Distributed Averaging}
\label{alg:global-average}
\end{algorithm}

\begin{lem}\label{lemma:convergence}
		\cite[Theorem 1]{Mou2014}: Suppose $\mathcal{G}$ is a tree graph with diameter equal to $d$. Algorithm \ref{alg:global-average} makes it possible for each agent $a$ to get the global {parameters} at a maximum of $d$ steps. 
	\end{lem}

The maximum path between two agents in the MDST obtained will be the maximum number of iterations in which the algorithm converges to global values ($P_D$, $P_G$, $N$ and $M$).  

\subsection{Distributed Projected Consensus Gradient} \label{subsec:DPG}

To solve the problem presented in \eqref{eqn:opt} we use the projected consensus algorithm proposed in \cite{Ozdaglar2010}. Our algorithm does not require private information from the neighboring agents such as its incremental cost. We only need the power estimates of generators and loads in iteration $k$ to estimate $k+1$. 

 \begin{assumption}
$P_{g_{i}}(0)$ and $P_{d_{j}}(0)$ for all {$i \in {G}$ and $j\in {D}$} are feasible points, i.e., $P_{g_{i}}(0)$ and $P_{d_{j}}(0)$ satisfy \eqref{eqn:gen_constraints} and \eqref{eqn:con_constraints}. 
\label{as:initial_values}
\end{assumption}   
   
For the initialization step we let Assumption \ref{as:op_set_ne}, \ref{as:graph}, \ref{as:neighbors}  and \ref{as:initial_values} hold. {The $a$-th agent updates its estimate by using the information produced by Algorithm \ref{alg:global-average},  then taking a gradient step to minimize the cost or maximize the utility function, and then projecting the result onto its constraint set $X_a$, where $\mathit{X_{a}}$ is the set of feasible solutions for each agent (cf. Eqs. \eqref{eqn:gen_constraints} and \eqref{eqn:con_constraints}).} 
Initially, we seek for reach an average consensus, for this we use the Algorithm \ref{alg:global-average}. We define the stacked vector of power (generated and demanded) as $P^\top=[P_{g}, P_{d}]^\top$, $P_a$ is the power of agent $a \in \mathcal{V}$, $P_b$ is the power of the neighbors of agent $a$, with $b\in \mathcal{N}_a$, we use the variable $v_a$ to store the sum of the powers in the iterative system as follows

\begin{align}
	\left\{ 
    \begin{array}{ll}
	v_a(1)= & P_{a}(0)+ \sum_{b\in \mathcal{N}_{a}}P_{b}(0)\\
	v_a(2)= & \sum_{b\in \mathcal{N}_{a}}P_{b}(1)+(1-r_a)P_{a}(0)\\
	 	    & \vdots					\\
    v_a(d)= & \sum_{b\in \mathcal{N}_{a}}P_{b}(d-1)+(1-r_a)P_{D_a}(d-2).\\
	\end{array} 
	\right.
    \label{eqn:sum_P}
\end{align}

When \eqref{eqn:sum_P} has been executed, every agent $a \in \mathcal{V}$ has a value  $v_a(d)=  \sum_{a\in \mathcal{V}}P_{a}$. Once $v_a(d)$ is obtained, it is possible to take the gradient step. For this step, we use $z(k)$ which contains the power average consensus minus the gradient of cost or utility function such as

\begin{align}
	z(k)=\frac{v_a(d)}{N_a+M_a}-\alpha_k d_a(k),
    \label{eqn:updatelaw}
\end{align}

where $\alpha_{k} > 0$ is the stepsize, $d_a(k)$ is the gradient of $U(P_{d_j}(k))$ and $C(P_{g_i}(k))$ depending on each agent. Finally, $z(k)$ is projected onto the feasible sets $X_a$,  The projection vector is denoted as $\mathbb{P}_{X_a}[\cdot]$ and it is defined as $\mathbb{P}_{\mathit{X_a}}[Y]= \argmin_{x \in \mathit{X_a}}  \Vert x-Y \Vert $. Each agent makes projections taking into account the constraints to which it is subjected. Consumers and generators are subject to constraints associated to its maximum and minimum load and generation, respectively. The projection onto the feasible set is defined as follows:

\begin{align}\label{eqn:projections}
	\begin{split}
	\mathbb{P}_{X_a}[z(k)]=z_1(k)= 
	\begin{cases}
	\overline{P_{a}},	& \text{if}\qquad  z(k)>\overline{P_{a}},\\
	\underline{P_{a}},	& \text{if}\qquad  z(k)<\underline{P_{a}},\\
	z(k) & \text{Otherwise.}
	\end{cases}
	\end{split}
\end{align}

Furthermore, generators have to maintain the global constraint shown in \eqref{eqn:global_cons}, the projection onto the feasible set $X_g$, where $X_{g}$ is the constraint set where constraint \eqref{eqn:global_cons} is held, is shown in \eqref{eqn:projection_generators}.

\begin{align}\label{eqn:projection_generators}
	\begin{split}
	\mathbb{P}_{X_g}[z_1(k)]= \hspace{5cm}\\
	\begin{cases}
	z_1(k),	& \text{if \eqref{eqn:global_cons} is hold},\\
	z_1(k)-\dfrac{P_{G_a}(k)-P_{D_a}(k)}{N_a} & \text{Otherwise}.
	\end{cases}
	\end{split}
\end{align}

Finally, the power update law $P_a(k+1)$ is given by for consumers and generator agents as follows

\begin{align}
\begin{split}
	P_{a}(k+1)=
    \begin{cases}
    \mathbb{P}_{X_a}[z(k)] & \text{if } a \in D \\
    \mathbb{P}_{X_g}[z_1(k)] &  \text{if } a \in G \label{eqn:p(k+1)_gen}
	\end{cases}
\end{split}
\end{align}

 \subsection{ Distributed Transactive Algorithm}	

We now state the distributed transactive algorithm, that is, the main contribution on this paper. We use $k$ to denote the iterations in the algorithm.

\begin{algorithm} 
 \begin{algorithmic}[1]
 	\State \textbf{Executed by:}  Agents $a \in \mathcal{V}= \{1,...,N+M\}$
 	\State \textbf{Require:} $d_a$, $\alpha_k$ and $r_a$
    \State \textbf{Initialize:} Assumption \ref{as:initial_values} is hold $\forall$ $a$.  Set $k=0$
    \State \textbf{Execute} Algorithm MDST in Subsection \ref{subsec:DAMDST}.
 		\While{$k\geq0$}
        \State \textbf{Send} $P_{a}(k)$ to all $b \in \mathcal{N}_a$             
 		\State \textbf{Execute} Algorithm \ref{alg:global-average}
        \State \textbf{Obtain} $P_G(k)$, $P_D(k)$, $M(k)$ and $N(k)$
 		\State \textbf{Execute} For all $a \in \mathcal{V}$\\
        				\hspace{3cm}Equation \eqref{eqn:sum_P} and then \eqref{eqn:updatelaw}
   		\State \textbf{Execute} For all $a \in \mathcal{V}$\\
        				\hspace{3cm}Equation \eqref{eqn:projections} 
        \State \textbf{Execute} For all $i \in G$\\
        				\hspace{3cm}Equation \eqref{eqn:projection_generators}
        \State \textbf{Obtain} $P_a(k+1)$
        \State\textbf{Set} $k=k+1$;
 		\EndWhile
 \end{algorithmic}
\caption{Distributed Transactive Control Algorithm}
\label{alg:distributed_transactive_control}
\end{algorithm}

\subsection{Convergence}\label{subsec:convergence}

In this subsection, we analyze and prove the convergence of the proposed distributed transactive algorithm. 

\begin{thm}\label{thm:optimum}
Assume that the stepsize $\alpha$ satisfies that $\sum_k \alpha_k = \infty $ and $\sum_k \alpha^{2}_{k} \geq \infty $. 
Furthermore, let ${P_a(k)}$, with $a \in \mathcal{V}$, be the set points generated by Algorithm \ref{alg:distributed_transactive_control} and $X = \cap^{N+M}_{a=1}X_a$ be the intersection set between all feasible sets of the agents. Then, $P_a(k)$ with $a \in \mathcal{V}$ converges to the optimal solution $P_{a}^{*}$ with $P_{a}^{*} \in X$, that is

\begin{align*}
		\lim_{k\to\infty} P_a(k) = P_{a}^{*}.
		\end{align*}
\end{thm}

\begin{pf} 
		Without loss of generality, all agents can be listed such as in \eqref{eqn:opt_generalized}. The optimization problem defined in \eqref{eqn:opt} can be generalized as

		\begin{mini!}[2]
			{P_{a}}{\sum_{a=1}^W -U_{a}(P_{a})} 	
			{\label{eqn:opt_generalized}}{}
			\addConstraint{\sum_{a=1}^W P_{a}}{=0} \label{eqn:global_cons_generalized}
			\addConstraint{\overline{P_{a}}\geq P_{a} \geq \underline{ P_{a}}}{\quad \forall a  \in W}\label{eqn:limits_generalization} 
		\end{mini!}

where $W=N+M$, and

\begin{subequations}
			\begin{align}
			P_a(k)= 
			\begin{cases}
			P_{d_j},	&  j=1,..., M\\
			-P_{g_{i+M}}, 	&  i=M+1,...,W
			\end{cases}
			\end{align}
			\begin{align}
			U_a(\cdot)= 
			\begin{cases}
			U_j(\cdot),		&  j=1,..., M\\
			-C_i(\cdot), 		&  i=M+1,...,W
			\end{cases}
			\end{align}
		\end{subequations}
	
 The constraints \eqref{eqn:Tc_gen_constraints} and \eqref{eqn:Tc_con_constraints} are satisfied through the stepsize $\alpha_k$. Let $X_{t_{a}}$ be the set where \eqref{eqn:Tc_gen_constraints} and \eqref{eqn:Tc_con_constraints} are feasible.  Let $\mathcal{D}$ be the diameter of the set $X_{t_{a}}$, i.e., $\mathcal{D}=\max_{x,y}||x-y||$, where $x, y\in X_{t_{a}}$. Then, there exists an $\alpha_k$ such that $\alpha_k d_a(k) \in \mathcal{D}\hspace{0.15cm} \forall \hspace{0.15cm} k$ and thus satisfying restrictions \eqref{eqn:Tc_gen_constraints} and \eqref{eqn:Tc_con_constraints}. Moreover, constraints \eqref{eqn:global_cons_generalized} and \eqref{eqn:limits_generalization} can be written such as $X = \cap^{W}_{a=1}X_a \cap X_g$, where $X_g$ is the set of feasible solutions for \eqref{eqn:global_cons_generalized}. It is possible to write \eqref{eqn:opt_generalized} such as

		\begin{mini!}[2]
			{P_{a}}{\sum_{a=1}^W -U_{a}(P_{a})} 	
			{\label{eqn:theorem_Nedic}}{}
			\addConstraint{X = \cap^{W}_{a=1}X_a \cap X_g}{} 
		\end{mini!}

It is assumed that $X_a$ and $X_g$ are compact sets. Given that all utility functions $U_{a}(P_{a})$ are continuous, and based on Weierstrass' Theorem \eqref{eqn:theorem_Nedic} has an optimal solution $P_{a}^{*} \in X$. { Considering that $\frac{v_a(d)}{N_a+M_a}$ in \eqref{eqn:updatelaw} is equal to $\lim_{k\to\infty}\frac{1}{N+M}\sum_{a\in M+N}P_a$ and} according to   \cite[Proposition 5]{Ozdaglar2010} agents executing Algorithm \ref{alg:distributed_transactive_control} converges to $P_{a}^{*} \in X$. $\blacksquare$
        \end{pf}

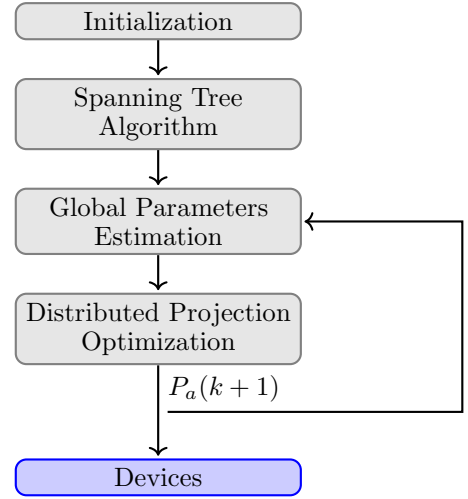
\begin{figure}[tb]
			\centering
            \begin{tikzpicture} [
    auto,
    block/.style    = { rectangle, draw=gray, thick, 
                        fill=gray!20, text width=10em, text centered,
                        rounded corners, minimum height=1em },
   device/.style    = { rectangle, draw=blue, thick, 
                        fill=blue!20, text width=10em, text centered,
                        rounded corners, minimum height=1em },
    line/.style     = { draw, thick, ->, shorten >=1pt },
  ]
  \matrix [column sep=4mm, row sep=5mm] {
                    & \node [block] (initialization) {Initialization};              & \\
                    & \node [block] (SPT) {Spanning Tree Algorithm};   & \\
                    & \node [block] (GPE) {Global Parameters Estimation}; & \\
                    & \node [block] (DPO){Distributed Projection Optimization};& \\
                    & \node (null1) {};                                    & \\
              		& \node [device] (other_agent){Devices}; \\
  };
  \begin{scope} [every path/.style=line]
    \path (initialization)        --    (SPT);
    \path (SPT)      			  --    (GPE);
    \path (null1)   --++  (4,0)  |- (GPE);
    \path (GPE)    		  --    (DPO);
   \path (DPO) 	      --    node [near start] {$P_a(k+1)$} (other_agent);
    \end{scope}
\end{tikzpicture}
			\caption{Flowchart of distributed transactive algorithm}
			\label{fig:c}
\end{figure}
    \section{Case Studies}\label{sec:cases}

In this section, we simulate a distribution system with five distributed generators and five consumers able to change its loads, depending on the system state. We seek to maximize the social welfare of the population of generator and consumers. Agents have limited power generation and demand, and consumers and generators satisfy Assumption 2. In order to simplify $\underline{T_{a}}$ and $\overline{T_{a}}$ for all agent, they are assumed as $-100$ and $100$ respectively.

	\begin{table}[htb]
		\caption{System Parameters in Simulation}
\begin{center}
		\label{table:Agents_data}
		\begin{tabular}{@{}ccc|ccc@{}}
\toprule
\multicolumn{1}{l}{\textbf{Generator}} & \multicolumn{1}{l}{$\mathbf{\overline{P_{g_i}}}$} & \multicolumn{1}{l|}{$\mathbf{\underline{P_{g_i}}}$} & \multicolumn{1}{l}{\textbf{Consumers}} & \multicolumn{1}{l}{$\mathbf{\overline{P_{d_j}}}$} & \multicolumn{1}{l}{$\mathbf{\underline{P_{d_j}}}$} \\ \midrule
1                     & 4000                & 100                      & 1                                      & 4100                                & 3000                                 \\
2                    & 6000               & 100                       & 2                                      & 5200                                & 4000                                 \\
3                   & 7000                 & 100                     & 3                                      & 6300                                & 5000                                \\
4                & 8000                 & 100                        & 4                                      & 6400                                & 5000                                \\
5                   & 9000                  & 100               & 5                                      & 7500                                & 6000                                \\
                  &                  &                            & 6                                      & 2000                                & 0                                \\ \bottomrule
\end{tabular}
	\end{center}
\end{table}

\subsection{Simulation with Smart Loads}

    \begin{figure*}[ht]
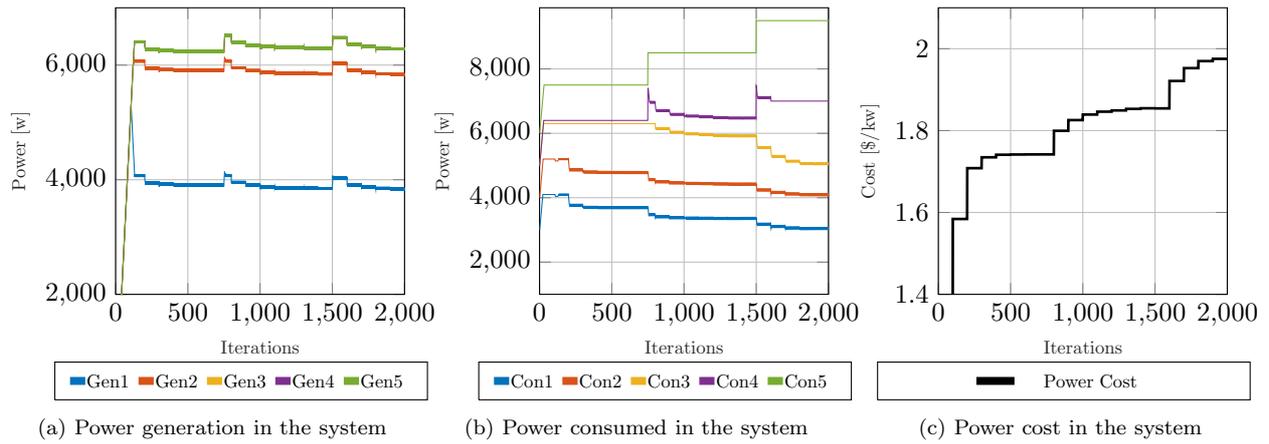

        \begin{subfigure}{.3\textwidth}
			\centering
            \input{generators.tex}
			\caption{Power generation in the system }
			\label{fig:case2_generation}
		\end{subfigure}
		\begin{subfigure}{.3\textwidth}
			\centering
            \input{loads.tex}
			\caption{Power consumed in the system}
			\label{fig:case2_demanded}
		\end{subfigure}
		\begin{subfigure}{.3\textwidth}
			\centering
            \input{globalcost.tex}            
			\caption{Power cost in the system}
			\label{fig:case2_cost}
		\end{subfigure}  
        \centering
		\caption{Simulations of system with Smart Loads}
		\label{fig:Case_2} 
	\end{figure*}

We use five generators and five smart loads and the parameters in Table \ref{table:Agents_data} for each agent. Smart loads and generators are connected indistinctly, i.e., it is not necessary that the system has a specific topology. Each $750$ iterations the base load changes in the simulation, it is considered the case where loads $4$ and $5$ rise up its base loads to $1000$ W. We refer to time instant $k$ as iterations. Figure \ref{fig:case2_generation}  shows the  power delivered by each generator, it is possible to see that instead of the changes in the load, generators can supply the exact power demanded. Figure \ref{fig:case2_demanded} shows the simulation results for the power demanded by Smart Loads. Figure \ref{fig:case2_demanded} shows that smart loads lower their consumption when the power cost in the system rises as a result of the increase in fixed load. The cost power is shown in Figure \ref{fig:case2_cost}. Figure \ref{fig:case2_cost} shows clearly that the price increases as more power is required by the loads, therefore adjustable loads reduce its consumption as is shown in Figure \ref{fig:case2_demanded}.

    \subsection{Adding a New Agent to the System: Smart Load}

\begin{figure*}[ht]
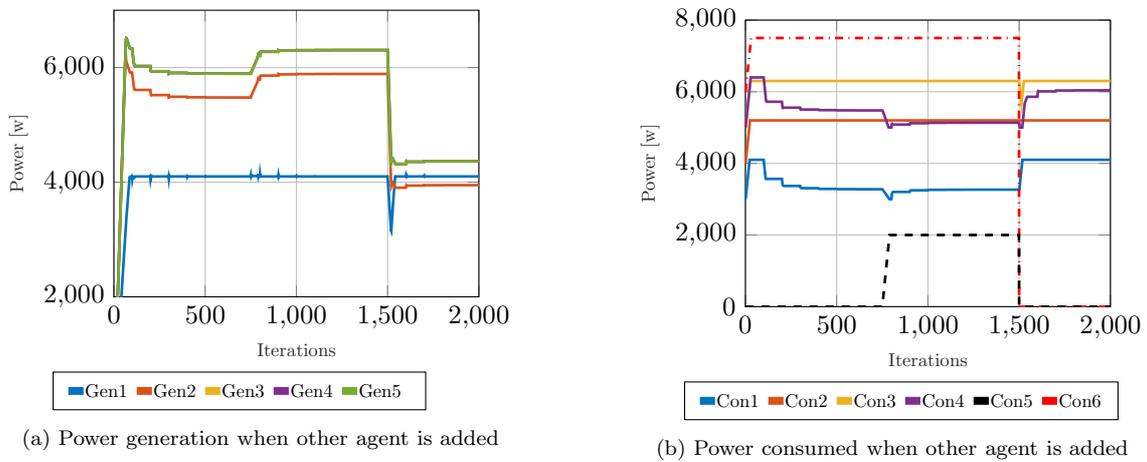

		\centering
		\begin{subfigure}{.45\textwidth}
			\centering
			\input{gen_1.tex}
			\caption{Power generation when other agent is added}
			\label{fig:case_P}
		\end{subfigure}
		\begin{subfigure}{.45\textwidth}
			\centering
            \input{loads_1.tex}
			\caption{Power consumed when other agent is added}
			\label{fig:case2_P_na}
		\end{subfigure} 
		\caption{System simulations when other agent consumer is added to the system}
		\label{fig:Case_32} 
	\end{figure*}

In this subsection we add a smart load to the system, the new agent is the consumer $6$. When we add this agent to the system, we assume that the system can recognize it and execute Algorithm MDST in subsection \ref{subsec:DAMDST}. It is possible observe in Figure \ref{fig:case_P} and Figure \ref{fig:case2_P_na} that we add the new agent at $750$ iterations since the new agent is configured in the system and demand power, the cost of the energy in the system will increase, it is possible to observe that generation matches exactly the power demanded.  We remove at $1500$ iterations two agents to the system, consumers $5$ and $6$, leaving the system with $5$ generators and $4$ consumers, for this reason, power cost decreases and variable loads of agent $1$, $2$ and $4$ rise its load. Despite the rise of load, the total load demanded by the system when two agents are removed is lower than the load in the previous state. Therefore, in Figure \ref{fig:case_P} power generated decreases matching the demanded power.
	
\section{Conclusion and Future Work}\label{sec:conclusions}

We proposed a new control strategy in the dispatch of distributed generators and demand of the users in power systems based on a transactive control framework. We consider some constraints in the generators and consumers with satisfactory results. Besides, we demonstrate that distributed transactional controllers are capable of addressing problems with distributed information in power systems. The simulation results show that the distributed transactional control algorithm achieves optimal social welfare in a dynamic way while maintaining system constraints in a power network. The study of adversarial agents in the power system for transactive control requires future study.

\bibliography{ifacconf}           
\end{document}